\documentstyle[amscd,amssymb,verbatim]{amsart}
\theoremstyle{plain}
\newtheorem{Thm}{Theorem}

\newtheorem{Lem}{Lemma}

\numberwithin{equation}{section}
\begin{document}
\title[Landen transformations]
{A geometric view of the rational Landen transformations\\}
\author{John Hubbard }
\address{Department of Mathematics, Cornell University,
Ithaca, New York, NY}
\email{jhh8@@cornell.edu}
\author{Victor Moll}
\address{Department of Mathematics, Tulane  University,
New Orleans, LA 70118}
\email{vhm@@math.tulane.edu}
\subjclass{Primary: 33; Secondary: 65}
\date{March 21, 2001}

\maketitle

\medskip

\begin{abstract}
We provide a geometric interpretation of a new rational 
Landen transformation and establish the convergence 
of its iterates.
\end{abstract} 

\medskip

\newcommand{\nn}{\nonumber}
\newcommand{\no}{\noindent}
\newcommand{\cl}{\mathop{\rm Cl}\nolimits}
\newcommand{\li}{\mathop{\rm Li}\nolimits}
\newcommand{\realpart}{\mathop{\rm Re}\nolimits}
\newcommand{\imagpart}{\mathop{\rm Im}\nolimits}
\newcommand{\lp}{\ln \sqrt{2 \pi}}
\newcommand{\ba}{\begin{eqnarray}}
\newcommand{\ea}{\end{eqnarray}}
\newcommand{\ift}{\int_{0}^{\infty}}
\newcommand{\ione}{\int_{0}^{1}}
\newcommand{\npR}{\mathbb{R}^{-}_{0}}
\newcommand{\nnR}{\mathbb{R}^{+}_{0}}
\newcommand{\allR}{\mathbb{R}}
\newcommand{\allN}{\mathbb{N}}
\newcommand{\nnN}{\mathbb{N}_{0}}
\newcommand{\an}{{\bf{a}}_{n}}
\newcommand{\bn}{{\bf{b}}_{n}}
\newcommand{\ano}{{\bf{a}}}
\newcommand{\bno}{{\bf{b}}}
\newcommand{\azero}{{\bf{a}}_{0}}
\newcommand{\bzero}{{\bf{b}}_{0}}
\newcommand{\pone}{{\mathbb{P}}^{1}}
\newcommand{\pz}{\pone{\kern-4pt_{z}}}
\newcommand{\pw}{\pone{\kern-4pt_{w}}}
\newcommand{\pd}{B^{-}{\kern-4pt_{d}} (\pone)}
\newcommand{\vf}{\varphi}

\section{Introduction}
\setcounter{equation}{0}

The transformation theory of elliptic integrals was initiated by Landen in
\cite{landen1,landen2}, wherein he proved the invariance of the function
\ba
G(a,b) & = & \int_{0}^{\pi/2} \frac{d \theta}{\sqrt{a^{2} \cos^{2} \theta +
b^{2} \sin^{2} \theta}} \label{elli1}
\ea
\no
under the transformation
\ba
a_{1} = \frac{a+b}{2} & \quad  & b_{1} = \sqrt{ab}. \label{ite1}
\ea
\no
Gauss \cite{gauss1} rediscovered this invariance in the process of
calculating the arclength of  a lemniscate. The limit of the sequence $(a_{n},
b_{n})$ defined by iteration of (\ref{ite1}) is the celebrated
{\em arithmetic-geometric mean} $\text{AGM}(a,b)$ of $a$ and $b$. The 
invariance of
the elliptic integral (\ref{elli1}) leads to
\ba
\frac{\pi}{2 \, \text{AGM}(a,b)} & = & G(a,b). \label{inv2}
\ea
\no
General information about AGM and its applications is given in 
\cite{borbor}. A geometric interpretation of the transformation (\ref{ite1}) 
is given in \cite{grayson}. \\

A transformation analogous to the Gauss-Landen map (\ref{ite1}) has been
given in \cite{bm2} for the rational integral
\ba
U_{6}(a_{1},a_{2};b_{0},b_{1},b_{2}) & = &
\ift \frac{b_{0}z^{4} + b_{1}z^{2} + b_{2} }
{z^{6} + a_{1}z^{4} + a_{2}z^{2} + 1} \, dz. \label{u6}
\ea
Indeed, the integral $U_{6}$ is invariant under the transformation
\ba
a_{1}^{(1)} & = & \frac{a_{1}a_{2} + 5a_{1}+ 5a_{2} + 
9}{(a_{1}+a_{2}+2)^{4/3}} \label{formula1} \\
a_{2}^{(2)} & = & \frac{a_{1}+a_{2}+6}{(a_{1}+a_{2}+2)^{2/3}} \nn \\
b_{0}^{(1)} & = & \frac{b_{0}+b_{1}+b_{2}}{(a_{1}+a_{2}+2)^{2/3}} \nn \\
b_{1}^{(1)} & = & \frac{b_{0}(a_{2}+2) + 2b_{1} + b_{2}(a_{1}+3)}
{a_{1} + a_{2} + 2} \nn \\
b_{2}^{(1)} &  = & \frac{b_{0} + b_{2}}{(a_{1}+a_{2}+2)^{2/3}}.  \nn
\ea
\medskip

This transformation was obtained by a sequence of elementary changes of
variables and the convergence of
\ba
(\an,\bn) & := & (a_{1}^{(n)},a_{2}^{(n)},b_{0}^{(n)},b_{1}^{(n)},b_{2}^{(n)}).
\nn
\ea
\no 
For any initial data $(\azero,\bzero)
\in \mathbb{R}_{+}^{2} \times \mathbb{R}_{+}^{3}$ there exists a number $L$,
depending upon the initial condition, such that
\ba
(\an, \bn) & \longrightarrow & (3,3,L,2L,L),  \label{conv1}
\ea
\no
so that
\ba
U_{6}(\an,\bn) & \longrightarrow & L \times \frac{\pi}{2}.
\ea
\no
The invariance of $U_{6}$ under (\ref{formula1}) shows that
\ba
U_{6}(\azero,\bzero) & = & L \times \frac{\pi}{2}. \label{conv2}
\ea
\no
Therefore the iteration given above becomes an iterative form of evaluating the
integral. \\

The main result of \cite{bm3}, quoted below,  is an extension of 
 (\ref{formula1}) for an even integrand. 

\begin{Thm}
\label{thmbig}
Let $R(z) = P(z)/Q(z)$ with
\ba
P(z) = \sum_{j=0}^{p-1} b_{j} z^{2(p-1-j)} \text{   and   } \; \;
Q(z) = \sum_{j=0}^{p} a_{j} z^{2(p-j)}. \label{polyb}
\ea
\no
Define $a_{j}=0$ for $j>p$, $b_{j}=0$ for $j>p-1$, 
\ba
d_{p+1-j} & = & 
\sum_{k=0}^{j} a_{p-k} a_{j-k} 
\label{dofj}
\ea
\no
for $0 \leq k \leq p-1$,
\ba
d_{1} & = & 
\frac{1}{2} \sum_{k=0}^{p} a_{p-k}^{2},
\label{dofj1}
\ea
\ba
c_{j} & = & \sum_{k=0}^{2p-1} a_{j}b_{p-1-j+k} 
\label{cofj}
\ea
\no
for $0 \leq j \leq  2p-1$, and
\ba
\alpha_{p}(i) & = & \begin{cases}
     2^{2i-1} \sum_{k=1}^{p+1-i} \frac{k+i-1}{i} \binom{k+2i-2}{k - 1}
d_{k + i} \text{    if } 1 \leq i \leq p    \\
     1 + \sum_{k=1}^{p} d_{k} \text{ if } i = 0. 
          \end{cases} 
\label{alpha}
\ea
\no
Let
\ba
a_{i}^{+} & = & \frac{\alpha_{p}(i)}{2^{2i} Q(1)^{2(1-i/p)}} 
\label{acoeff}
\ea
\no
for $1 \leq i \leq p-1$, and
\ba
b_{i}^{+} & = & Q(1)^{2i/p+1/p-2}  \times 
\left[ \sum_{k=0}^{p-1-i} (c_{k} + c_{2p-1-k}) \binom{p-1-k+i}{2i} 
\right] \label{bcoeff}  
\ea
\no
for $0 \leq i \leq p-1$. 
\no
Finally, define the polynomials 
\ba
P^{+}(z)  = {\sum_{k=0}^{p-1} b_{i}^{+} z^{2(p-1-i)}} & \text{ and } &
Q^{+}(z) = {\sum_{k=0}^{p} a_{i}^{+} z^{2(p-i)}}. \label{newb}
\ea
\no
Then 
\ba
\int_{0}^{\infty} \frac{P(z)}{Q(z)} \; dz & = & 
\int_{0}^{\infty} \frac{P^{+}(z)}{Q^{+}(z)} \; dz. 
\ea
\end{Thm}

\medskip

The proofs in \cite{bm2,bm3}  are elementary but lack a proper geometric
interpretation. In particular, the proof of (\ref{conv1}) given in \cite{bm2} 
could not be extended even for degree 8 in view of the formidable
algebraic difficulties. The goal of this paper is to show that the 
transformation (\ref{acoeff}, \ref{bcoeff})  is a
particular case of a general construction: the direct image of a meromorphic
$1$-form under a rational map. This will allow us to prove an analogue of
(\ref{conv1},\ref{conv2}) for  the  integral
\ba
U_{2p}(\ano,\bno) & := & \ift
\frac{b_{0}z^{2p-2} +b_{1}z^{2p-4} + \cdots + b_{p}}
{z^{2p} + a_{1}z^{2p-2} + \cdots + 1} \, dz. \label{up}
\ea
\no
In fact, we prove that the sequence ${\mathbf{x}}_{n}$ starting at 
\ba
{\mathbf{x}}_{0} & =  & (a_{1}, \cdots, a_{p-1}; b_{0}, \cdots, b_{p-1})  \nn
\ea
\no
and defined by ${\mathbf{x}}_{n+1} = {\mathbf{x}}_{n}^{+}$ satisfies 
\ba
{\mathbf{x}}_{n} & \to & \left( 
\binom{p}{1}, \binom{p}{2}, \cdots, \binom{p}{p-1}; 
\binom{p-1}{0} \, L, \binom{p-1}{1} \, L, \cdots, 
\binom{p-1}{p-1} \, L \right),  \nn
\ea
\no
where
\ba
L & = & \frac{2}{\pi} \, U_{2p}(\ano,\bno).  \nn
\ea
\no
Moreover the convergence of the iteration is equivalent to 
the convergence of the initial integral. \\

\bigskip

\section{The direct image of a 1-form}
\setcounter{equation}{0}

Let $\pi: X \rightarrow Y$ be a proper analytic mapping of Riemann surfaces
(i.e., a finite ramified covering space),  and
  $\vf$ be a tensor of any type on $X$. Then $\pi_{*} \vf$ is the tensor of the
same  type on $Y$, defined as follows: Let $U \subset Y$ be a simply connected
subset of $Y$ containing no critical value of $\pi$, and let
$\sigma_{1}, \cdots, \sigma_{k}:U \rightarrow X$ be the distinct sections of
$\pi$. Then the 
direct image of $\pi_*\vf$ is defined by
\ba
\pi_{*} \vf\  \Big{|}_{U} &  = & \sum_{i=1}^{k} \sigma_{i}^{*} \vf.
\ea
\no
This defines $\pi_{*}\vf$ except at the ramification values of
$\pi$, where $\pi_{*}\vf$ may acquire poles even if $\vf$ is holomorphic.

We will be applying this construction in the case where $\vf$ is a holomorphic
1-form, and in this case $\pi_*\vf$ is analytic.

\begin{Lem}
\label{lemdirimagehol}
If $\pi:X \to Y$ is proper and analytic as above, and $\vf$ is an
analytic 1-form on $X$, then $\pi_*\vf$ is an analytic 1-form on $Y$.
Furthermore, for any 
oriented rectifiable curve $\gamma$ on $Y$, we have
$$
\int_\gamma \pi_* \vf = \int_{\pi^{-1}\gamma} \vf.
$$
\end{Lem}

\begin{pf} The only problem is to show that $\pi_*\vf$ is 
holomorphic at the critical
values.  It is clearly enough to show that the contribution of a 
neighborhood of a single
critical point is holomorphic.  Thus we may assume that 
$\pi(z)=w=z^m$ for some $m$,
and that
$$
\vf = (a_kz^k+ a_{k+1}z^{k+1}+ \dots) dz,
$$
with $k\ge 0$.

For $i=0, \dots, m-1$ set $\sigma_i(w) = \zeta^i \sigma_0(w)$, where 
$\zeta= e^{2\pi
i/m}$ and $\sigma_0(w) = w^{1/m}$ for some branch of the $1/m$ power, 
for instance
the one where the argument is between $0$ and $2\pi/m$.  Then
\ba
\pi_{*}(z^{k} dz)  & = &
\begin{cases}
0 \quad \hspace{50pt}  \text{ if } k+1 \text{ is not divisible by $m$} \\
u^{(k+1-m)/m} du  \quad \hspace{9pt}  \text{ if $k+1$ is divisible by $m$}.
\end{cases}
\no \label {compdirimagezm}
\ea

Thus the first term of the power series for $\vf$ to contribute 
anything to $\pi_*\vf$ is
the term of degree $m-1$, and it contributes to the constant term; similarly, 
the terms of
degree $2m-1, 3m-1, \dots$ contribute to the terms of degree $1,2, 
\dots$, all positive
powers.
\end{pf}

This has a useful corollary. Recall that the degree of a meromorphic 
function is the
maximum of the degrees of the numerator and the denominator when the 
rational function is
written in reduced form.

\begin{Lem}
\label{degsdirimratfuncs}
If $\pi:\pone \to \pone$ is  analytic, and $\vf = R(z) dz$ is a
meromorphic 1-form on $\pone$ so that $R$ is a rational function of 
degree $k$, then
$\pi_*\vf$ can be written as 
$R_1(z) dz$, where $R_1$ is a rational function of degree at most $k$.
\end{Lem}

\begin{pf}
By Lemma \ref{lemdirimagehol}, the number of poles of $\pi_*\vf$ is at 
most equal to the
number of poles of $\vf$, and clearly the orders of the poles cannot 
increase either.
\end{pf}

\noindent
{\bf Note. } It is quite possible for the degree of $\pi_*\vf$ to be 
less than the degree
of $\vf$.  This can happen in two ways: we might have poles at two 
points $z_1, z_2$ such
that $\pi(z_1)=\pi(z_2)$, and then the polar parts at these points 
could cancel.  We may
also have a pole of order $>1$ at a critical point, and then the 
order of the pole at the
corresponding critical value will be less (in fact, the pole might 
disappear altogether).

\section{A particular branched cover}
\setcounter{equation}{0}

We will be concerned with the specific map
\ba
\pi(z) & = & w := \frac{z^{2}-1}{2z}. \label{branch1}
\ea
\no

This mapping can also be viewed as the Newton map associated to the 
equation $z^2+1=0$.  As
such it has $\pm i$ as superattractive fixed points, and $\pi$
is conjugate to $F(z) = z^{2}$ via the Mobius transformation
$M(z) = (z+i)/(z-i)$; indeed $M \circ \pi \circ M^{-1} = F$.

Let us list some properties of $\pi$.

\begin {Lem} If
$\vf$ has no poles on
$\overline
\allR
\subset \pone$, then
$$
\int_{-\infty}^\infty \vf =\int_{-\infty}^\infty \pi_* \vf.
$$
\end{Lem}

\begin{pf} If $\vf$ has no poles on $\allR$ (including at 
infinity), then the integral
converges. Since $\pi$ maps the real axis (including $\infty$) to 
itself as a double
cover, the result follows from Lemma
\ref{lemdirimagehol}.
\end{pf}

Let $\tau:\pone \to \pone$ be the map $z \mapsto -z$.  Then clearly 
$\pi\circ \tau = \tau
\circ \pi$.  Call $\vf$ even if $\tau^*\vf = \vf$, and odd if 
$\tau^*\vf = -\vf$. \\

\no
{\bf Note}. When $\vf= R(z) \, dz$ with $R$ a rational function, then 
$\vf$ is even if and
only if $R$ is odd, and $\vf$ is odd if and only if $R$ is even, 
since $dz$ is odd.

\begin{Lem} We have the following identities:

(a) $$\pi^*\pi_* \vf = \vf + \tau^* \vf.
$$

(b) If $\vf$ is even, then $\pi_* \vf=0$.

(c) If $\vf$ is odd, then $\pi_* \vf $ is also odd.
\end{Lem}

Thus we can restrict our attention to odd 1-forms. Below we calculate 
$\pi_*(R(z) \, dz),$
where $R(z)$ is an even rational function.  We will only consider the 
case when the
numerator of $R$ has degree at least 2 less than the denominator, as 
this avoids a pole at
infinity which would prevent the integral over $\allR$ from converging.

The explicit evaluations of the form
$\pi_{*}\vf$ described below were conducted using Mathematica. The 
corresponding sections are 
\ba
\sigma_{\pm}(w)  & = & w \pm \sqrt{w^{2} + 1},
\ea
\no
so that for $\varphi = \Phi(z) dz$ we have 
\ba
\pi_{*} \varphi & = & \Phi(\sigma_{+}(w)) \frac{d \sigma_{+}}{dw} + 
\Phi(\sigma_{-}(w)) \frac{d \sigma_{-}}{dw}. 
\ea
\no
The calculations are formidable. \\

\no
{\bf Example 1}. Let
\ba
\varphi & = & \frac{b_{0}}{a_{0}z^{2} + a_{1}} \, dz.
\ea
\no
Then

\ba
\pi_{*}\varphi & = & \frac{2b_{0}(a_{0} + a_{1})}{4a_{0}a_{1}w^{2} +
(a_{0}+a_{1})^{2}} \, dw.
\ea

\no
Observe that the new $1$-form can be written as
\ba
\pi_{*} \varphi & = & b_{0} \times
\frac{A(a_{0},a_{1})}{G^{2}(a_{0},a_{1})w^{2} + A^{2}(a_{0},a_{1})} \, dw,
\ea

\no
where $A(a,b)$ and $G(a,b)$ are the arithmetic and geometric means of $a$ and
$b$ respectively.   \\

\bigskip

\no
{\bf Example 2}. The form
\ba
\varphi & = & \frac{b_{0}z^{2} + b_{1} }{a_{0}z^{4} + a_{1}z^{2} + a_{2}} \,
dz
\ea
\no
is transformed into \\

\ba
\pi_{*}\varphi & = &
\frac{
8(a_{2}b_{0} + a_{0}b_{1})w^{2} + 2(a_{0} + a_{1} + a_{2})(b_{0}+b_{1})}
{
16a_{0}a_{2}w^{4} + 4(a_{0}a_{1} + 4a_{0}a_{2} + a_{1}a_{2})w^{2} +
(a_{0}+a_{1}+a_{2})^{2} } \; dw. \nn
\ea

\bigskip

\section{The convergence of $(\pi_*)^n \vf$}
\setcounter{equation}{0}

In this section we present the principal theorem of the paper.

\begin{Thm}
\label{covergencethm}
   Let $\vf$ be a 1-form, holomorphic on a neighborhood $U$ of $\allR
\subset \pone$.  Then 
$$
\lim_{n \to \infty} (\pi_*)^n \vf= \frac 1\pi 
\left(\int_{-\infty}^\infty \vf\right)
\frac {dz}{1+z^2},
$$
\no
where the convergence is uniform on compact subsets of $U$.
\end{Thm}

\begin{pf} We find it convenient to prove this for the map 
$F(z)=z^2$, which is
conjugate to $\pi$.  In that form, the statement to be proved is that 
if $\vf$ is analytic
in some neighborhood $U$ of the unit circle, then
$$
\lim_{n \to \infty} (F_*)^n \vf = \frac 1{2\pi i} \left(\int_{S^1} 
\vf\right) \frac{dz}z.
$$

\no
Any such 1-form $\vf$ can be developed in a Laurent series
$$
\vf = \left(\sum_{k=-\infty}^\infty a_kz^k\right) \frac {dz}z,
$$
where $\sum_{k=1}^\infty (|a_k|+|a_{-k}|) \rho^k<\infty$ for some 
$\rho>1$.  Note that
$$
a_0 = \frac 1{2\pi i} \int_{S^1} \vf.
$$

\no
In this form it is very easy to compute $F_* \vf$.

\begin{Lem}
The mapping $F_*$ on 1-forms is given by
$$
F_* \vf = \sum_{k=-\infty}^\infty a_{2k} z^k \frac {dz}z.
$$
\end{Lem}

\begin{pf} This is what was computed in Equation \ref {compdirimagezm}.
\end{pf}

Thus in the ``basis'' of forms $z^k \frac {dz}z$, the vector 
corresponding to $k=0$ is
an eigenvector with eigenvalue $1$, and the rest of the space is nilpotent:
$$
(F_*)^m z^k \frac {dz}z=0
$$
if $m$ is greater than the greatest power of $2$ which divides $k$. 
This comes close to
proving Theorem \ref{covergencethm}, but it doesn't quite; for instance
$$
\left(\sum_{k=0}^\infty z^k\right) \frac {dz}z= \frac{dz}{z(1-z)}
$$
is also fixed under $F_*$.  We cannot argue merely in terms of formal 
Laurent series:
convergence must be taken into account.

But this is not too hard. Consider the region $U_R$ defined by
$ \frac 1 R <  |z| < R$,
and the space $A_R$ of analytic 1-forms
$$
\vf = \left(\sum_{k=-\infty}^\infty a_kz^k\right) \frac {dz}z
$$
on $U_R$ such that
$$
\|\phi\| =|a_0|+ \sum_{k=1}^\infty (|a_k|+a_{-k}|) R^{k}  < \infty .
$$

\no
We then have
\ba 
\|\pi_*^n\vf-a_0\frac {dz}z\| & = & \sum_{k=1}^\infty (|a_{2^nk}|+ 
|a_{-2^nk}|) R^{k} \nn \\
& =  & 
\sum_{k=1}^\infty (|a_{2^nk}|+ |a_{-2^nk}|) R^{2^nk}\frac {R^k}{R^{2^n}k} \nn \\
& \le & 
\frac R{R^{2^n}}\|\vf\|. \nn
\ea
This certainly shows that $\pi_*^n\vf-a_0\frac {dz}z$ tends to $0$, 
in fact very fast: it
superconverges to $0$.

\end{pf}

\bigskip

\section{Normalization of the integrands}
\setcounter{equation}{0}

In the previous section we have produced a map $\pi_*$ of 1-forms 
$\phi = R(z) \, dz$ that
does not increase the degree and the integral over $[0,\infty]$. 
Moreover, we have seen
that the integrands $\pi_*^n \vf$ converge as $n$ tends to infinity. 
This doesn't quite
tell us about the convergence of the coefficients of $R$, because of 
possible common
factors and cancellations.  Here we normalize the rational functions 
so that $\pi_*$
induces a convergent iteration on the coefficients.

We will write the  integrands so that their denominators are monic and
with constant term  equal to $1$. The latter can be achieved by 
factoring the constant
term out  while the former is obtained by a change of
variable of the
form $z \mapsto \lambda z$, with an appropriate $\lambda$. \\

\no
{\bf Example 1}. For rational functions of degree $2$ we obtain  \\

\ba
\ift \frac{b_{0}}{a_{0}z^{2} + a_{1}} \, dz & = &
\ift \frac{2b_{0}(a_{0}+a_{1}) }{4a_{0}a_{1}w^{2} + (a_{0}+a_{1})^{2} } \, dw.
\ea
\no
This is an identity: both sides normalize to
\ba
\frac{b_{0}}{\sqrt{a_{0}a_{1}}} \times \ift \frac{dx}{x^{2}+1}.
\ea

\bigskip

\no
{\bf Example 2}. The quartic case yields \\

\ba
\ift \frac{b_{0}z^{2} + b_{1} }{a_{0}z^{4} + a_{1}z^{2} + a_{2} } \; dz & = &
\ift \frac{b_{0}^{(1)}w^{2} + b_{1}^{(1)} }{a_{0}^{(1)}w^{4} + 
a_{1}^{(1)}w^{2} + a_{2}^{(1)} } \; dw,
\ea
\no
where
\ba
b_{0}^{(1)} & = & 8(a_{2}b_{0} + a_{0}b_{1}) \label{quartic1} \\
b_{1}^{(1)} & = & 2(a_{0} + a_{1} + a_{2})(b_{0} + b_{1}) \nn \\
a_{0}^{(1)} & = & 16a_{0}a_{2} \nn \\
a_{1}^{(1)} & = & 4(a_{0}a_{1} + 4a_{0}a_{2} + a_{1}a_{2}) \nn  \\
a_{2}^{(1)} & = & (a_{0}+a_{1}+a_{2})^{2}. \nn
\ea

\medskip

\no
The normalization shows that the integral

$$ \ift \frac{b_{0}a_{2}^{1/2} z^{2} + b_{1}a_{0}^{1/2}}
{z^{4} + a_{0}^{-1/2} a_{1}a_{2}^{-1/2} z^{2} + 1} \, dz $$

\no
equals

$$ (a_{0} + a_{1} + a_{2})^{-1/2} \times $$

$$ \ift \frac{(a_{2}b_{0} + a_{0}b_{1})w^{2} +
(b_{0}+b_{1})a_{0}^{1/2}a_{2}^{1/2} }
{ w^{4} + [a_{0}a_{1} + 4a_{0}a_{2} +
a_{1}a_{2})a_{0}^{-1/2}a_{2}^{-1/2}(a_{0}+a_{1}+a_{2})^{-1} ] w^{2} + 1}
\, dw. $$

\medskip

\no
Naturally this identity can be verified directly using 
$$
\ift \frac{dx}{x^4 + 2ax^2  + 1}  = 
\ift \frac{x^{2} \, dx}{x^4 + 2ax^2  + 1}  =  \frac{\pi}{2^{3/2} \sqrt{a+1}}.
$$

\bigskip

\no
{\bf Example 3}. In the case of degree $6$ we obtain  \\

\bigskip

\ba
\ift \frac{b_{0}z^{4} + b_{1}z^{2} + b_{2}}
{a_{0}z^{6} + a_{1}z^{4} + a_{2}z^{2} + a_{3}} \; dz & =  &
\ift \frac{b_{0}^{(1)}w^{4} + b_{1}^{(1)}w^{2} + b_{2}^{(1)}}
{a_{0}^{(1)}w^{6} + a_{1}^{(1)}w^{4} + a_{2}^{(1)}w^{2} + a_{3}^{(1)}} \; dw,
\nn \\
\label{landen6}
\ea
\no
where
\ba
b_{0}^{(1)} & = & 32(a_{3}b_{0} + a_{0}b_{2}) \label{degreesix} \\
b_{1}^{(1)} & = &  8(a_{2}b_{0} + 3a_{3}b_{0} + a_{0}b_{1} + a_{3}b_{1} +
3a_{0}b_{2} + a_{1}b_{2}) \nn \\
b_{2}^{(1)} & = & 2(a_{0} + a_{1}+a_{2} + a_{3})(b_{0} + b_{1} + b_{2}) \nn \\
a_{0}^{(1)} & = & 64a_{0}a_{3} \nn \\
a_{1}^{(1)} & = & 16(a_{0}a_{2} + 6a_{0} a_{3} + a_{1}a_{3}) \nn \\
a_{2}^{(1)} & = & 4(a_{0}a_{1} + 4a_{0}a_{2} + a_{1}a_{2} + 9a_{0}a_{3} +
4a_{1}a_{3} + a_{2}a_{3}) \nn \\
a_{3}^{(1)} & = & (a_{0} + a_{1} + a_{2} + a_{3})^{2}. \nn
\ea
\no
The normalization of (\ref{landen6}) yields (\ref{formula1}).  \\

\bigskip

\no
{\bf Acknowledgments}. The second author acknowledges the partial support 
of NSF-DMS 0070567, Project number 540623. This work originated on a 
field trip to the rain forest in Puerto Rico 
during SIMU 2000. The authors thank Herbert
Medina and Ivelisse Rubio for their hospitality.  \\

\end{document}